\documentclass[12pt]{amsart}
 
\pdfoutput=1
\usepackage[margin=1in]{geometry}
\usepackage{amsmath,amsthm,amssymb,amsrefs,bbm,color,esvect,float,graphicx}
\newenvironment{theorem}[2][Theorem]{\begin{trivlist}
\item[\hskip \labelsep {\bfseries #1}\hskip \labelsep {\bfseries #2.}]}{\end{trivlist}}
\newenvironment{lemma}[2][Lemma]{\begin{trivlist}
\item[\hskip \labelsep {\bfseries #1}\hskip \labelsep {\bfseries #2.}]}{\end{trivlist}}

\theoremstyle{definition}

\begin{document}

\title[Weak-Type Estimates for Multilinear Calder\'on-Zygmund Operators]{An Endpoint Weak-Type Estimate for Multilinear Calder\'on-Zygmund Operators}
\author{Cody B. Stockdale}
\address{Cody B. Stockdale, Department of Mathematics, Washington University in St. Louis, One Brookings Drive, St. Louis, MO, 63130, USA}
\email{codystockdale@wustl.edu}
\author{Brett D. Wick}
\address{Brett D. Wick, Department of Mathematics, Washington University in St. Louis, One Brookings Drive, St. Louis, MO, 63130, USA}
\email{wick@math.wustl.edu}
\thanks{B. D. Wick's research is supported in part by National Science Foundation grant DMS \#1560955 and \#1800057.}

\begin{abstract}
The purpose of this article is to provide an alternative proof of the weak-type $\left(1,\ldots,1;\frac{1}{m}\right)$ estimate for $m$-multilinear Calder\'on-Zygmund operators on $\mathbb{R}^n$ first proved by Grafakos and Torres. Subsequent proofs in the bilinear setting have been given by Maldonado and Naibo and also by P\'erez and Torres. The proof given here is motivated by the proof of the weak-type $(1,1)$ estimate for Calder\'on-Zygmund operators in the nonhomogeneous setting by Nazarov, Treil, and Volberg.
\end{abstract}
\maketitle
\noindent \textbf{Keywords:} singular integrals; multilinear operators; weak-type estimates.

\noindent \noindent MSC Primary 42B20

\section{Introduction}
The Calder\'on-Zygmund theory of singular integral operators is central in the study of harmonic analysis. A key property of Calder\'on-Zygmund operators in the linear setting is their boundedness from $L^p(\mathbb{R}^n)$ to $L^p(\mathbb{R}^n)$ for all $p\in(1,\infty)$, assuming a priori that the operators are bounded from $L^2(\mathbb{R}^n)$ to $L^2(\mathbb{R}^n)$. The well-known method of proof is as follows:
\begin{enumerate}
\item establish a weak-type $(1,1)$ estimate for the operator, 
\item use the Marcinkiewicz interpolation theorem to obtain a strong type $(p,p)$ bound for all $p \in (1,2)$, and 
\item use duality to deduce the strong type $(p,p)$ estimate for all $p \in (1,\infty)$.
\end{enumerate}
For a detailed treatment of the classical Calder\'on-Zygmund theory, see \cites{Grafakos,Stein}.

The classical proof of the weak-type $(1,1)$ bound utilizes the Calder\'on-Zygmund decomposition for $L^1(\mathbb{R}^n)$ functions. This technique readily extends to handle Calder\'on-Zygmund operators on spaces that have an underlying measure possessing the doubling property. Such spaces are called spaces of homogeneous type. Recall that a measure $\mu$ possesses the doubling property if there exists a constant $C>0$ such that $$\mu(B(x,2r))\leq C\mu(B(x,r))$$ for all $r>0$ and all $x$ in the space. The classical technique, however, does not generalize as easily to spaces of nonhomogeneous type, which are spaces whose underlying measures instead have a polynomial growth condition. To address this setting, Tolsa developed a version of the Calder\'on-Zygmund decomposition adapted to nonhomogeneous measures to prove the weak-type estimate in a similar manner to the classical case in \cite{Tolsa2001}. In \cite{NTV1998}, Nazarov, Treil, and Volberg provided a proof of the weak-type $(1,1)$ bound of Calder\'on-Zygmund operators in the nonhomogeneous setting without using the Calder\'on-Zygmund decomposition. The proof in \cite{NTV1998} also works in the classical setting on $\mathbb{R}^n$.

More recently, attention has been given to the study of multilinear Calder\'on-Zygmund operators (see \cites{DamianLernerPerez2015,GrafakosKalton2001,GrafakosTorres2002,LOPTT-G2009,Lin2016,MN2009,PerezTorres2014}). To describe the setting, let $m$ be a positive integer. We say $K:\mathbb{R}^{n(m+1)} \rightarrow \mathbb{C}$ is an \emph{$m$-multilinear Calder\'on-Zygmund kernel} if there exist $C_K, \delta >0$ such that the following conditions hold:
\begin{enumerate}
\item (\emph{size}) \[|K(x,y_1,\ldots,y_m)| \leq \frac{C_K}{\left(\sum_{i=1}^m|x-y_i|\right)^{nm}}\]
for all $x,y_1,\ldots,y_m \in \mathbb{R}^n$ with $x\neq y_j$ for some $j$,
\item (\emph{smoothness}) \[|K(x,y_1,\ldots,y_m)-K(x',y_1,\ldots,y_m)| \leq \frac{C_K|x-x'|^{\delta}}{(\sum_{i=1}^m|x-y_i|)^{nm+\delta}}\] whenever $\displaystyle |x-x'| \leq \frac{1}{2}\max_{1\leq i\leq m}|x-y_i|$, and 
\[|K(x,y_1,\ldots,y_j,\ldots ,y_m)-K(x,y_1,\ldots, y_j', \ldots ,y_m)|\leq \frac{C_K|y_j-y_j'|^{\delta}}{(\sum_{i=1}^m|x-y_i|)^{nm+\delta}}\] for each $j\in \{1,\ldots,m\}$ whenever $\displaystyle |y_j-y_j'| \leq \frac{1}{2}\max_{1\leq i\leq m}|x-y_i|$.
\end{enumerate}
We say a bounded multilinear operator $T:(L^2(\mathbb{R}^n))^m \rightarrow L^{\frac{2}{m}}(\mathbb{R}^n)$ is an \emph{$m$-multilinear Calder\'on-Zygmund operator associated to a kernel $K$} if $K$ is an $m$-multilinear Calder\'on-Zygmund kernel and if \[T(f_1,\ldots,f_m)(x) = \int_{(\mathbb{R}^n)^m} K(x,y_1,\ldots,y_m)f_1(y_1)\cdots f_m(y_m)dy_1\cdots dy_m\] for almost every $x \in \mathbb{R}^n \setminus \bigcap_{i=1}^m\text{supp }f_i$. Throughout the remainder of this paper, $T$ will denote a multilinear Calder\'on-Zygmund operator. 

Let $\mathcal{M}(\mathbb{R}^n)$ denote the space of $\mathbb{R}$-valued Borel measures on $\mathbb{R}^n$. For $\nu_1,\ldots, \nu_m \in \mathcal{M}(\mathbb{R}^n)$ and $x \in \mathbb{R}^n \setminus \bigcap_{i=1}^m\text{supp }\nu_i$, define \[T(\nu_1,\ldots, \nu_m)(x):=\int_{(\mathbb{R}^n)^m}K(x,y_1,\ldots, y_m)d\nu(y_1)\cdots d\nu(y_m).\] We will denote the total variation of $\nu \in \mathcal{M}(\mathbb{R}^n)$ by $\|\nu\|$. Notice that if $f$ is a Borel measurable function, then $\||f|dm\|=\|f\|_{L^1(\mathbb{R}^n)}$ and \[T(f_1dm,\ldots,f_mdm)=T(f_1,\ldots,f_m).\] Here $fdm\in\mathcal{M}(\mathbb{R}^n)$ is defined for Borel subsets of $\mathbb{R}^n$ by $fdm(A)=\int_A f(x)dx$. Also, if $\nu_{i}=\sum_{j=1}^{N}a_{i,j}\delta_{x_{i,j}}$ for $i\in\{1,\ldots,m\}$, then \[T(\nu_1,\ldots,\nu_m)(x)=\sum_{j_1,\ldots,j_m=1}^{N}\left(\prod_{i=1}^ma_{i,j_i}\right)K(x,x_{1,j_1},\ldots,x_{m,j_m}).\]

Analogous properties to the classical case were established for multilinear operators in \cite{GrafakosTorres2002} (see also \cites{MN2009,PerezTorres2014}). In particular, a weak-type estimate is proved and used to establish strong type estimates via interpolation. The appropriate weak-type estimate for multilinear Calder\'on-Zygmund operators is of type $\left(1,\ldots,1;\frac{1}{m}\right)$. It is stated as the following: 
\begin{theorem}{1}
Let $T$ be a multilinear Calder\'on-Zygmund operator. If $f_1, \ldots, f_m \in L^1(\mathbb{R}^n)$, then \[\|T(f_1,\ldots,f_m)\|_{L^{\frac{1}{m},\infty}(\mathbb{R}^n)} \leq A_2\prod_{i=1}^m\|f_i\|_{L^1(\mathbb{R}^n)}\] where $A_2$ depends on $C_K$, $n$, and $m$.
\end{theorem}
As in the classical situation, Grafakos and Torres \cite{GrafakosTorres2002} prove Theorem 1 using the Calder\'on-Zygmund decomposition. An alternative proof is presented in Section 3.

The new proof is modeled after the argument for the weak-type estimate in \cite{NTV1998}. Instead of obtaining cancellation by means of the Calder\'on-Zygmund decomposition, we do so by subtracting terms involving certain point mass measures. The argument is then completed by establishing a weak-type estimate on a mixture of linear combinations of point mass measures and of $L^1(\mathbb{R}^n)\cap L^{\infty}(\mathbb{R}^n)$ functions with appropriate $L^{\infty}(\mathbb{R}^n)$ norm. This is stated precisely as the following: 
\begin{theorem}{2}
Let $T$ be a multilinear Calder\'on-Zygmund operator, $t>0$, and $l \in \{1,\ldots,m\}$ be given. If $\nu_1,\ldots, \nu_l \in \mathcal{M}(\mathbb{R}^n)$ are of the form $\nu_i = \sum_{j=1}^{N} a_{i,j}\delta_{x_{i,j}}$ where $x_{i,j} \in \mathbb{R}^n$ and $a_{i,j}\in\mathbb{R}$ and if $f_{l+1},\ldots,f_{m}\in L^1(\mathbb{R}^n)\cap L^{\infty}(\mathbb{R}^n)$ satisfy $\|f_i\|_{L^{\infty}(\mathbb{R}^n)}\leq t^{\frac{1}{m}}$ for all $i$, then \[|\{|T(\nu_1,\ldots,\nu_{l},f_{l+1},\ldots,f_{m})|>t\}|\leq A_3t^{-\frac{1}{m}}\left(\prod_{i=1}^l\|\nu_i\|^{\frac{1}{m}}\right)\left(\prod_{i=l+1}^{m}\|f_i\|_{L^1(\mathbb{R}^n)}^{\frac{1}{m}}\right)\] where $A_3$ depends on $C_K$, $n$, and $m$.
\end{theorem} It is not important that the $\nu_i$ are applied in the first $l$ slots of $T$ -- an identical proof yields the theorem whenever the set of indices of the $\nu_i$ is a nonempty subset of $\{1,\ldots,m\}$.

Since the proof of Theorem 1 still requires a decomposition of the arbitrary $L^1(\mathbb{R}^n)$ functions into bounded and unbounded pieces (``good'' and ``bad'' pieces), there is an analogy to be made between our proof and the proof in \cite{GrafakosTorres2002}. First, the term where the operator is only being applied to the ``good'' functions is handled identically -- both proofs use Chebyshev's inequality, the a priori boundedness of $T$, and the $L^{\infty}(\mathbb{R}^n)$ norms of the good functions to obtain the appropriate estimate. However, the terms where the operator has at least one ``bad'' function as an input are treated differently.

First we describe the Calder\'on-Zygmund decomposition approach to handling these terms. Because of the nature of this decomposition, each ``bad'' function, $b_i$, can be written as the sum $b_i=\sum_{j=1}^{\infty}b_{i,j}$ where each $b_{i,j}$ has mean value zero, is supported on a cube of appropriate measure, and has useful $L^1(\mathbb{R}^n)$ control. The cancellation involved in the $b_{i,j}$ allows one to introduce a term with the kernel evaluated at the center of the cube on which $b_{i,j}$ is supported, then one can use the smoothness assumption of the kernel to obtain the desired estimate. The disjointness of $\text{supp}(b_{i,j})$ over the $j$ allows one to recover the estimate for the original term.

Without the Calder\'on-Zygmund decomposition, there is no immediate cancellation that may be exploited in the ``bad'' functions. Instead, we apply a Whitney decomposition to write the support of each ``bad'' function, $b_i$, as a union of dyadic cubes with disjoint interiors and restrict $b_i$ to each cube given by the Whitney decomposition. Call these restrictions $b_{i,j}$. It suffices to approximate $b_i \approx \sum_{j=1}^{N} b_{i,j}$ and get an appropriate estimate with these sums, uniform in $N$. With this goal in mind, denote the center of $\text{supp}(b_{i,j})$ by $c_{i,j}$ and define measures $\nu_i^N$ by $\nu_i^N:=\sum_{j=1}^N a_{i,j}\delta_{c_{i,j}}$ where $a_{i,j} := \int_{\mathbb{R}^n}b_{i,j}(x) dx$. Adding and subtracting terms involving these $\nu_i^N$ introduces cancellation. We then subtract a term involving the kernel evaluated at the $c_{i,j}$ and use the regularity discussed in Section 2 to get the desired control. It is then left to control a term involving a mixture of linear combinations of point mass measures and ``good'' functions. This can be done with Theorem 2.

Section 2 includes Lemma 1, a regularity condition first described in \cite{PerezTorres2014} for bilinear kernels, which is the key use of cancellation. Section 3 contains the main results. The proof of Theorem 1 assuming Theorem 2 is given first. The proof of Theorem 2 is given at the end.

We would like to acknowledge Rodolfo Torres for kindly providing comments and references.


\section{Preliminaries}
The proofs of Theorem 1 and Theorem 2 use the multilinear geometric H\"ormander condition given in Lemma 1 below. This type of regularity was first introduced in the bilinear setting by P\'erez and Torres in \cite{PerezTorres2014}. Throughout the rest of this paper, we use the notation $\vv{y}_{i,k}=(y_i,y_{i+1},\ldots,y_k)$, $\vv{f}_{i,k}=(f_i,f_{i+1},\ldots,f_k)$, $\vv{\nu}_{i,k}=(\nu_i,\nu_{i+1},\ldots,\nu_{k})$, $\vv{c}_{(i,j_i),(k,j_k)}=(c_{i,j_i},c_{i+1,j_{i+1}},\ldots,c_{k,j_k})$, and $\vv{\nu}_{(i,j_i),(k,j_{k})}=(\nu_{i,j_i},\nu_{i+1,j_{i+1}},\ldots,\nu_{k,j_k})$. We apologize for further complicating the notation; however, this is necessary to compactly describe many expressions that follow.
\begin{lemma}{1}
\label{lemma1} There exists $A_1>0$ such that if $l\in\{1,\ldots,m\}$ and $\mathcal{S}_1,\ldots,\mathcal{S}_l$ are countable collections of sets satisfying either 
\begin{enumerate}
\item each $\mathcal{S}_i=\{S_{i,1},S_{i,2},\ldots\}$ consists of dyadic cubes with disjoint interiors or
\item each $\mathcal{S}_i=\{S_{i,1},S_{i,2},\ldots\}$ consists of sets satisfying: 
\begin{itemize}
\item $S_{i,j}$ have disjoint interiors,
\item $S_{i,j}\subseteq B(c_{i,j},r_{i,j})$, and 
\item $\displaystyle \Omega_i :=\bigcup_{j=1}^{\infty} S_{i,j}=\bigcup_{j=1}^{\infty} B(c_{i,j},r_{i,j})$,
\end{itemize}
\end{enumerate}
then
\begin{align*}
&\sum_{j_1,\ldots,j_l=1}^{\infty}\prod_{i=1}^l|S_{i,j_i}|\\
&\quad\quad\times\int_{\mathbb{R}^{n(m-l)}}\sup_{\substack{(y_1,\ldots,y_l) \\\in \prod_{i=1}^lS_{i,j_i}}} \int_{\mathbb{R}^n\setminus\left(\bigcup_{i=1}^l \Omega_{i}^{*}\right)} |K(x,\vv{y}_{1,m})-K(x,\vv{c}_{(1,j_1),(l,j_l)},\vv{y}_{l+1,m})|dxd\vv{y}_{l+1,m}\\
&\leq A_1\sum_{i=1}^{l}\left|\Omega_{i}\right|
\end{align*} 
where $\Omega_i^{*}:= \bigcup_{j=1}^{\infty} 2S_{i,j}=\bigcup_{j=1}^{\infty}B(c_{i,j},2r_{i,j})$.
\end{lemma} 
It is not important that the indices of the $\mathcal{S}_i$ range from $1$ to $l$ -- an identical proof yields the lemma whenever the set of indices is a nonempty subset of $\{1,\ldots,m\}$.

This regularity was considered in \cite {PerezTorres2014} when the $\mathcal{S}_i$ are collections of dyadic cubes with disjoint interiors. We will use the lemma when the collections $\mathcal{S}_i$ consist of dyadic cubes in the proof of Theorem 1 and when they are of the second type in the proof of Theorem 2.

\begin{proof} We only prove the statement when the collections $\mathcal{S}_i$ are of the second type. The proof for collections of dyadic cubes is similar and is addressed in the bilinear setting in \cite{PerezTorres2014}. For $i=1,\ldots, l$, fix $S_{i,j_i}\in \mathcal{S}_i$. Use the smoothness condition of $K$ and the fact that $S_{i,j_i} \subseteq B(c_{i,j_i},r_{i,j_i})\subseteq \overline{B(c_{i,j_i},r_{i,j_i})}$ to see
\begin{align*} 
\sup_{\substack{(y_1,\ldots,y_l) \\\in \prod_{i=1}^lS_{i,j_i}}} &\int_{\mathbb{R}^n\setminus\left(\bigcup_{i=1}^l \Omega_{i}^{*}\right)}|K(x,\vv{y}_{1,m})-K(x,\vv{c}_{(1,j_1),(l,j_l)},\vv{y}_{l+1,m})|dx\\
&\lesssim \sup_{\substack{(y_1,\ldots,y_l) \\\in \prod_{i=1}^lS_{i,j_i}}} \int_{\mathbb{R}^n\setminus\left(\bigcup_{i=1}^l \Omega_{i}^{*}\right)} \frac{\sum_{i=1}^l|y_i-c_{i,j_i}|^{\delta}}{(\sum_{i=1}^m|x-y_i|)^{nm+\delta}}dx\\
&\leq \sup_{\substack{(y_1,\ldots,y_l) \\\in \prod_{i=1}^l\overline{B(c_{i,j_i},r_{i,j_i})}}} \int_{\mathbb{R}^n\setminus\left(\bigcup_{i=1}^l \Omega_{i}^{*}\right)} \frac{\sum_{i=1}^lr_{i,j_i}^{\delta}}{(\sum_{i=1}^m|x-y_i|)^{nm+\delta}}dx.
\end{align*}
Since for fixed $y_i \in\overline{B(c_{i,j_i},r_{i,j_i})}$, $i=l+1,\ldots,m$, the function $\displaystyle\int_{\mathbb{R}^n\setminus\left(\bigcup_{i=1}^l \Omega_{i}^{*}\right)} \frac{\sum_{i=1}^lr_{i,j_i}^{\delta}}{(\sum_{i=1}^m|x-y_i|)^{nm+\delta}}dx$ is continuous in the variables $y_i\in \overline{B(c_{i,j_i},r_{i,j_i})}$, $i=1,\ldots, l$, and since $\overline{B(c_{i,j_i},r_{i,j_i})}$ is a compact set, we may write 
\begin{align*}
\sup_{\substack{(y_1,\ldots,y_l) \\\in \prod_{i=1}^l\overline{B(c_{i,j_i},r_{i,j_i})}}} &\int_{\mathbb{R}^n\setminus\left(\bigcup_{i=1}^l \Omega_{i}^{*}\right)} \frac{\sum_{i=1}^lr_{i,j_i}^{\delta}}{(\sum_{i=1}^m|x-y_i|)^{nm+\delta}}dx\\
&=\int_{\mathbb{R}^n\setminus\left(\bigcup_{i=1}^l \Omega_{i}^{*}\right)} \frac{\sum_{i=1}^lr_{i,j_i}^{\delta}}{(\sum_{i=1}^l|x-y_i^*|+\sum_{i=l+1}^m|x-y_i|)^{nm+\delta}}dx \quad\text{and}
\end{align*}
\begin{align*}
\inf_{\substack{(y_1,\ldots,y_l) \\\in \prod_{i=1}^l\overline{B(c_{i,j_i},r_{i,j_i})}}} &\int_{\mathbb{R}^n\setminus\left(\bigcup_{i=1}^l \Omega_{i}^{*}\right)} \frac{\sum_{i=1}^lr_{i,j_i}^{\delta}}{(\sum_{i=1}^m|x-y_{i}|)^{nm+\delta}}dx\\
&=\int_{\mathbb{R}^n\setminus\left(\bigcup_{i=1}^l \Omega_{i}^{*}\right)} \frac{\sum_{i=1}^lr_{i,j_i}^{\delta}}{(\sum_{i=1}^l|x-y_{i_*}|+\sum_{i=l+1}^m|x-y_i|)^{nm+\delta}}dx.
\end{align*}
Note that for $x\in\mathbb{R}^n\setminus\left(\bigcup_{i=1}^l \Omega_{i}^{*}\right)$, $|x-y_{i_*}|\leq 2r_i+|x-y_i^*|$ and $|x-y_i^*|\ge r_i$, so \[\frac{\sum_{i=1}^l|x-y_{i_*}|+\sum_{i=l+1}^m|x-y_i|}{\sum_{i=1}^l|x-y_i^*|+\sum_{i=l+1}^m|x-y_i|}\leq \frac{\sum_{i=1}^l2r_i}{\sum_{i=1}^l|x-y_i^*|+\sum_{i=l+1}^m|x-y_i|}+1\leq 3.\]
Then
\begin{align*}
&\sup_{\substack{(y_1,\ldots,y_l) \\\in \prod_{i=1}^lS_{i,j_i}}} \int_{\mathbb{R}^n\setminus\left(\bigcup_{i=1}^l \Omega_{i}^{*}\right)}|K(x,\vv{y}_{1,m})-K(x,\vv{c}_{(1,j_1),(l,j_l)},\vv{y}_{l+1,m})|dx\\
&\quad\lesssim \int_{\mathbb{R}^n\setminus\left(\bigcup_{i=1}^l \Omega_{i}^{*}\right)} \frac{\sum_{i=1}^lr_{i,j_i}^{\delta}}{(\sum_{i=1}^l|x-y_{i_*}|+\sum_{i=l+1}^m|x-y_{i}|)^{nm+\delta}}\\
&\quad\quad\quad\times\left(\frac{\sum_{i=1}^l|x-y_{i_*}|+\sum_{i=l+1}^m|x-y_i|}{\sum_{i=1}^l|x-y_i^*|+\sum_{i=l+1}^m|x-y_i|}\right)^{nm+\delta}dx\\
&\quad\lesssim \int_{\mathbb{R}^n\setminus\left(\bigcup_{i=1}^l \Omega_{i}^{*}\right)} \frac{\sum_{i=1}^lr_{i,j_i}^{\delta}}{(\sum_{i=1}^l|x-y_{i_*}|+\sum_{i=l+1}^m|x-y_i|)^{nm+\delta}}dx\\
&\quad=\inf_{\substack{(y_1,\ldots,y_l) \\\in \prod_{i=1}^l\overline{B(c_{i,j_i},r_{i,j_i})}}} \int_{\mathbb{R}^n\setminus\left(\bigcup_{i=1}^l \Omega_{i}^{*}\right)} \frac{\sum_{i=1}^lr_{i,j_i}^{\delta}}{(\sum_{i=1}^m|x-y_{i}|)^{nm+\delta}}dx\\
&\quad\leq \inf_{\substack{(y_1,\ldots,y_l) \\\in \prod_{i=1}^lS_{i,j_i}}} \int_{\mathbb{R}^n\setminus\left(\bigcup_{i=1}^l \Omega_{i}^{*}\right)} \frac{\sum_{i=1}^lr_{i,j_i}^{\delta}}{(\sum_{i=1}^m|x-y_{i}|)^{nm+\delta}}dx.
\end{align*}
Using the previous estimate, trivial estimates to pull the infimum outside of the integral, Fubini, and integral estimates, we get the bound
\begin{align*}
\int_{\mathbb{R}^{n(m-l)}}&\sup_{\substack{(y_1,\ldots,y_l) \\\in \prod_{i=1}^lS_{i,j_i}}} \int_{\mathbb{R}^n\setminus\left(\bigcup_{i=1}^l \Omega_{i}^{*}\right)}|K(x,\vv{y}_{1,m})-K(x,\vv{c}_{(1,j_1),(l,j_l)},\vv{y}_{l+1,m})|dxd\vv{y}_{l+1,m}\\
&\lesssim \int_{\mathbb{R}^{n(m-l)}}\inf_{\substack{(y_1,\ldots,y_l) \\\in \prod_{i=1}^lS_{i,j_i}}} \int_{\mathbb{R}^n\setminus\left(\bigcup_{i=1}^l \Omega_{i}^{*}\right)} \frac{\sum_{i=1}^lr_{i,j_i}^{\delta}}{(\sum_{i=1}^m|x-y_{i}|)^{nm+\delta}}dxd\vv{y}_{l+1,m}\\
&\leq \inf_{\substack{(y_1,\ldots,y_l) \\\in \prod_{i=1}^lS_{i,j_i}}} \int_{\mathbb{R}^n\setminus\left(\bigcup_{i=1}^l \Omega_{i}^{*}\right)}\int_{\mathbb{R}^{n(m-l)}} \frac{\sum_{i=1}^lr_{i,j_i}^{\delta}}{(\sum_{i=1}^m|x-y_i|)^{nm+\delta}}d\vv{y}_{l+1,m}dx\\
&\lesssim \inf_{\substack{(y_1,\ldots,y_l) \\\in \prod_{i=1}^lS_{i,j_i}}} \int_{\mathbb{R}^n\setminus\left(\bigcup_{i=1}^l \Omega_{i}^{*}\right)}\frac{\sum_{i=1}^lr_{i,j_i}^{\delta}}{(\sum_{i=1}^l|x-y_i|)^{nl+\delta}}dx.
\end{align*}
Therefore
\begin{align*}
&\sum_{j_1,\ldots,j_l=1}^{\infty}\prod_{i=1}^l|S_{i,j_i}|\\
&\quad\quad\times\int_{\mathbb{R}^{n(m-l)}}\sup_{\substack{(y_1,\ldots,y_l) \\\in \prod_{i=1}^lS_{i,j_i}}} \int_{\mathbb{R}^n\setminus\left(\bigcup_{i=1}^l \Omega_{i}^{*}\right)} |K(x,\vv{y}_{1,m})-K(x,\vv{c}_{(1,j_1),(l,j_l)},\vv{y}_{l+1,m})|dxd\vv{y}_{l+1,m}\\
&\lesssim C_K\sum_{j_1,\ldots,j_l=1}^{\infty}\prod_{i=1}^l|S_{i,j_i}|\inf_{\substack{(y_1,\ldots,y_l) \\\in \prod_{i=1}^lS_{i,j_i}}} \int_{\mathbb{R}^n\setminus\left(\bigcup_{i=1}^l \Omega_{i}^{*}\right)}\frac{\sum_{i=1}^lr_{i,j_i}^{\delta}}{(\sum_{i=1}^l|x-y_i|)^{nl+\delta}}dx\\
&\leq C_K \sum_{j_1,\ldots,j_l=1}^{\infty}\int_{S_{l,j_l}}\cdots\int_{S_{1,j_1}}\int_{\mathbb{R}^n\setminus\left(\bigcup_{i=1}^l \Omega_{i}^{*}\right)}\frac{\sum_{i=1}^lr_{i,j_i}^{\delta}}{(\sum_{i=1}^l|x-y_i|)^{nl+\delta}}dxd\vv{y}_{1,l}\\
&=C_K\sum_{k=1}^l\left(\sum_{\substack{j_1,\ldots,j_l=1 \\ r_{k,j_k}\ge r_{i,j_i}\,\text{all }i}}^{\infty}\int_{S_{l,j_l}}\cdots\int_{S_{1,j_1}}\int_{\mathbb{R}^n\setminus\left(\bigcup_{i=1}^l \Omega_{i}^{*}\right)}\frac{\sum_{i=1}^lr_{i,j_i}^{\delta}}{(\sum_{i=1}^l|x-y_i|)^{nl+\delta}}dxd\vv{y}_{1,l}\right).
\end{align*}

We will control the term of the summation above with $k=1$; the other terms are handled identically. Using trivial estimates, Fubini's theorem, the fact that the $S_{i,j_i}$ have disjoint interiors, and integral estimates, we obtain
\begin{align*}
&\sum_{\substack{j_1,\ldots,j_l=1 \\ r_{1,j_1}\ge r_{i,j_i}\,\text{all }i}}^{\infty}\int_{S_{l,j_l}}\cdots\int_{S_{1,j_1}}\int_{\mathbb{R}^n\setminus\left(\bigcup_{i=1}^l \Omega_{i}^{*}\right)}\frac{\sum_{i=1}^lr_{i,j_i}^{\delta}}{(\sum_{i=1}^l|x-y_i|)^{nl+\delta}}dxd\vv{y}_{1,l}\\
&\quad\lesssim \sum_{j_1,\ldots,j_l=1}^{\infty}\int_{S_{l,j_l}}\cdots\int_{S_{1,j_1}}\int_{\mathbb{R}^n\setminus\left(\bigcup_{i=1}^l \Omega_{i}^{*}\right)}\frac{r_{1,j_1}^{\delta}}{(\sum_{i=1}^l|x-y_i|)^{nl+\delta}}dxd\vv{y}_{1,l}\\
&\quad\lesssim \sum_{j_1=1}^{\infty}\int_{S_{1,j_1}}\int_{\mathbb{R}^n\setminus\left(\bigcup_{i=1}^l \Omega_{i}^{*}\right)}\sum_{j_2,\ldots,j_l=1}^{\infty}\int_{S_{l,j_l}}\cdots\int_{S_{2,j_2}}\frac{1}{(\sum_{i=1}^l|x-y_i|)^{nl}}d\vv{y}_{2,l}\frac{r_{1,j_1}^{\delta}}{|x-c_{1,j_1}|^{\delta}}dxdy_1\\
&\quad\leq \int_{\Omega_1}\int_{\mathbb{R}^n\setminus\left(\bigcup_{i=1}^l \Omega_{i}^{*}\right)}\int_{\mathbb{R}^{n(l-1)}}\frac{1}{(\sum_{i=1}^l|x-y_i|)^{nl}}d\vv{y}_{2,l}\frac{r_{1,j_1}^{\delta}}{|x-c_{1,j_1}|^{\delta}}dxdy_1\\
&\quad\lesssim \int_{\Omega_1}\int_{\mathbb{R}^n\setminus\left(\bigcup_{i=1}^l \Omega_{i}^{*}\right)}\frac{r_{1,j_1}^{\delta}}{|x-y_1|^n|x-c_{1,j_1}|^{\delta}}dxdy_1\\
&\quad\lesssim |\Omega_1|\int_{\mathbb{R}^n\setminus\Omega_{1}^{*}}\frac{r_{1,j_1}^{\delta}}{|x-c_{1,j_1}|^{n+\delta}}dx\\
&\quad\leq |\Omega_1|\int_{|x|>2r_{1,j_1}}\frac{r_{1,j_1}^{\delta}}{|x|^{n+\delta}}dx\\
&\quad\leq |\Omega_1|.
\end{align*}
Similarly, for $k=2,\ldots,l$, \[\displaystyle\sum_{\substack{j_1,\ldots,j_l=1 \\ r_{k,j_k}\ge r_{i,j_i}\,\text{all }i}}^{\infty}\int_{S_{l,j_l}}\cdots\int_{S_{1,j_1}}\int_{\mathbb{R}^n\setminus\left(\bigcup_{i=1}^l \Omega_{i}^{*}\right)}\frac{\sum_{i=1}^lr_{i,j_i}^{\delta}}{(\sum_{i=1}^l|x-y_i|)^{nl+\delta}}dxd\vv{y}_{1,l} \lesssim |\Omega_k|.\] This completes the proof.
\end{proof}


\section{Main Results}
We now turn to proving the main result of Grafakos and Torres \cite{GrafakosTorres2002}.
\begin{theorem}{1}
\label{theorem1}
If $f_1, \ldots, f_m \in L^1(\mathbb{R}^n)$, then \[\|T(f_1,\ldots,f_m)\|_{L^{\frac{1}{m},\infty}(\mathbb{R}^n)} \leq A_2\prod_{i=1}^m\|f_i\|_{L^1(\mathbb{R}^n)}\] where $A_2$ depends on $K$, $n$, $m$, and $\|T\|_{(L^2(\mathbb{R}^n))^m\rightarrow L^{\frac{2}{m}}(\mathbb{R}^n)}$. That is, for every $t>0$, it holds that \[|\{|T(f_1,\ldots,f_m)|>t\}|\leq A_2^{\frac{1}{m}}t^{-\frac{1}{m}}\prod_{i=1}^m\|f_i\|_{L^1(\mathbb{R}^n)}^{\frac{1}{m}}.\]
\end{theorem}
Our contribution is the following.
\begin{theorem}{2}
Let $t>0$ and $l \in \{1,\ldots,m\}$ be given. If $\nu_1,\ldots, \nu_l \in \mathcal{M}(\mathbb{R}^n)$ are of the form $\nu_i = \sum_{j=1}^{N} a_{i,j}\delta_{x_{i,j}}$ where $x_{i,j} \in \mathbb{R}^n$ and $a_{i,j}\in\mathbb{R}$ and if $f_{l+1},\ldots,f_{m}\in L^1(\mathbb{R}^n)\cap L^{\infty}(\mathbb{R}^n)$ satisfy $\|f_i\|_{L^{\infty}(\mathbb{R}^n)}\leq t^{\frac{1}{m}}$ for all $i$, then \[|\{|T(\nu_1,\ldots,\nu_{l},f_{l+1},\ldots,f_{m})|>t\}|\leq A_3t^{-\frac{1}{m}}\left(\prod_{i=1}^l\|\nu_i\|^{\frac{1}{m}}\right)\left(\prod_{i=l+1}^{m}\|f_i\|_{L^1(\mathbb{R}^n)}^{\frac{1}{m}}\right)\] where $A_3$ depends on $K$, $n$, $m$, and $\|T\|_{(L^2(\mathbb{R}^n))^m\rightarrow L^{\frac{2}{m}}(\mathbb{R}^n)}$.
\end{theorem}
Note that Theorem 2 holds whenever the set of indices of the $\nu_i$ is a nonempty subset of $\{1,\ldots,m\}$. We will first prove Theorem 1 assuming Theorem 2. We will then prove Theorem 2. Let $M$ denote the uncentered Hardy-Littlewood maximal function and recall its formula \[M f(x) := \sup_{Q \ni x} \frac{1}{|Q|} \int_{Q}|f(y)|dy.\] Write $\|M\|$ for $\|M\|_{L^1(\mathbb{R}^n)\rightarrow L^{1,\infty}(\mathbb{R}^n)}$.

\begin{proof}[Proof of Theorem 1]
Let $t>0$ be given. By density, we may assume $f_1,\ldots, f_m$ are continuous functions with compact support. Normalize to assume $\|f_1\|_{L^1(\mathbb{R}^n)}=\cdots=\|f_m\|_{L^1(\mathbb{R}^n)}=1$. Set \[G_i:=\left\{Mf_i>t^{\frac{1}{m}}\right\}\,\,\,\,\,\, \text{and}\,\,\,\,\,\, G:=\bigcup_{i=1}^mG_i.\] Notice that 
\[|G|\leq \sum_{i=1}^m |G_i| \leq m\|M\|t^{-\frac{1}{m}}.\] Put \[b_i:=f_i\mathbbm{1}_{G_i} \quad\quad \text{and}\quad\quad g_i:=f_i\mathbbm{1}_{\mathbb{R}^n\setminus G_i}.\] Set
\begin{align*}
E_1&:=\left\{\left|T\left(g_1,g_2,\ldots,g_m\right)\right|>\frac{t}{2^m}\right\},\\
E_2&:=\left\{\left|T\left(b_1,g_2,\ldots,g_m\right)\right|>\frac{t}{2^m}\right\},\\
E_3&:=\left\{\left|T\left(g_1,b_2,\ldots,g_m\right)\right|>\frac{t}{2^m}\right\},\\
&\quad\quad\quad\quad\quad\quad\quad\quad\cdots\\
E_{2^m}&:=\left\{\left|T\left(b_1,b_2,\ldots,b_m\right)\right|>\frac{t}{2^m}\right\};
\end{align*}
where each $E_s=\left\{\left|T(h_1,\ldots,h_m)\right|>\frac{t}{2^m}\right\}$ with $h_i\in\left\{b_i,g_i\right\}$ and all the sets $E_s$ are distinct. Since \[\left|\left\{\left|T(f_1,\ldots,f_m)\right|>t\right\}\right|\leq \sum_{s=1}^{2^m}|E_s|,\] it suffices to control each $|E_s|$ by a constant multiplied by $t^{-\frac{1}{m}}$.

We will first estimate $|E_1|$. Note that since $|f_i(x)|\leq Mf_i(x)\leq t^{\frac{1}{m}}$ for almost every $x \in \mathbb{R}^n\setminus G_i$, it is true that $\|g_i\|_{L^{\infty}(\mathbb{R}^n)}\leq t^{\frac{1}{m}}$. Use Chebyshev's inequality, the boundedness of $T$ from $(L^2(\mathbb{R}^n))^m$ to $L^{\frac{2}{m}}(\mathbb{R}^n)$, and the fact that $\|g_i\|_{L^{\infty}(\mathbb{R}^n)}\leq t^{\frac{1}{m}}$ to see
\begin{align*}
\left|E_1\right|&\leq 4t^{-\frac{2}{m}}\int_{\mathbb{R}^n}|T(g_1,\ldots,g_m)(x)|^{\frac{2}{m}}dx\\
&\leq 4\|T\|^{\frac{2}{m}}_{(L^2(\mathbb{R}^n))^m\rightarrow L^{\frac{2}{m}}(\mathbb{R}^n)}t^{-\frac{2}{m}}\prod_{i=1}^m\left(\int_{\mathbb{R}^n}|g_i(x)|^2dx\right)^{\frac{1}{m}}\\
&\leq 4\|T\|^{\frac{2}{m}}_{(L^2(\mathbb{R}^n))^m\rightarrow L^{\frac{2}{m}}(\mathbb{R}^n)}t^{-\frac{1}{m}}\prod_{i=1}^m\|g_i\|_{L^1(\mathbb{R}^n)}^{\frac{1}{m}}\\
&\leq 4\|T\|^{\frac{2}{m}}_{(L^2(\mathbb{R}^n))^m\rightarrow L^{\frac{2}{m}}(\mathbb{R}^n)}t^{-\frac{1}{m}}\prod_{i=1}^m\|f_{i}\|_{L^1(\mathbb{R}^n)}^{\frac{1}{m}}\\
&=4\|T\|^{\frac{2}{m}}_{(L^2(\mathbb{R}^n))^m\rightarrow L^{\frac{2}{m}}(\mathbb{R}^n)}t^{-\frac{1}{m}}\\
&=B_1t^{-\frac{1}{m}}
\end{align*}
where $B_1:=4\|T\|^{\frac{2}{m}}_{(L^2(\mathbb{R}^n))^m\rightarrow L^{\frac{2}{m}}(\mathbb{R}^n)}$.

Consider the set $E_s$ for a fixed $2\leq s\leq 2^m$. Suppose that there are $l$ functions of the form $b_i$ and $m-l$ functions of the form $g_i$ appearing as entries in the $T(h_1,\ldots,h_m)$ involved in the definition of $E_s$. For notational simplicity, assume that the $b_i$ are in the first $l$ entries and the $g_i$ are in the remaining $m-l$ entries (analogous arguments hold in the other cases). Apply a Whitney decomposition to write each $G_i$ as a union of dyadic cubes with disjoint interiors, \[G_i=\bigcup_{j=1}^{\infty}Q_{i,j},\] where \[2\text{diam}(Q_{i,j})\leq d(Q_{i,j},\mathbb{R}^n\setminus G_i) \leq 8\text{diam}(Q_{i,j}).\] Put \[b_{i,j}=b_i\mathbbm{1}_{Q_{i,j}} \,\,\,\, \text{and} \,\,\,\,  b_i^{N}=\sum_{j=1}^Nb_{i,j}.\]  It suffices to control (uniformly in $N$) the measure of $E_s$ with $b_i$ replaced by $b_i^{N}$. Denote this set by $\widetilde{E}_s$.

Let $c_{i,j}$ denote the center of $Q_{i,j}$ and set \[a_{i,j}=\int_{Q_{i,j}}b_{i,j}(x)dx,\quad\quad \nu_{i,j}=a_{i,j}\delta_{c_{i,j}}, \quad\quad \text{and} \quad\quad \nu_i^N=\sum_{j=1}^N \nu_{i,j}.\] 
Notice, by adding and subtracting $T\left(\vv{\nu^N}_{1,k},\vv{b^N}_{k+1,l},\vv{g}_{l+1,m}\right)$ for $1\leq k \leq l$, we have
\begin{align*}
\left|\widetilde{E}_s\right|&\leq\sum_{k=1}^{l}\left|\left\{\left|T\left(\vv{\nu^N}_{1,k-1},b_{k}^N-\nu_k^N,\vv{b^N}_{k+1,l},\vv{g}_{l+1,m}\right)\right|>\frac{t}{(l+1)2^m}\right\}\right|\\
&\quad\quad+\left|\left\{\left|T\left(\vv{\nu^N}_{1,l},\vv{g}_{l+1,m}\right)\right|>\frac{t}{(l+1)2^m}\right\}\right|\\
&\leq \sum_{k=1}^{l}\left(|G|+\left|\left\{\mathbb{R}^n\setminus G: \left|T\left(\vv{\nu^N}_{1,k-1},b_{k}^N-\nu_{k}^N,\vv{b^N}_{k+1,l},\vv{g}_{l+1,m}\right)\right|>\frac{t}{(l+1)2^m}\right\}\right|\right)\\
&\quad\quad+\left|\left\{\left|T\left(\vv{\nu^N}_{1,l},\vv{g}_{l+1,m}\right)\right|>\frac{t}{(l+1)2^m}\right\}\right|\\
&\leq m^2\|M\|t^{-\frac{1}{m}}+\sum_{k=1}^l|S_k|+|S|,
\end{align*}
where
\begin{align*}
&S_k:=\left\{\mathbb{R}^n\setminus G: \left|T\left(\vv{\nu^N}_{1,k-1},b_{k}^N-\nu_{k}^N,\vv{b^N}_{k+1,l},\vv{g}_{l+1,m}\right)\right|>\frac{t}{(l+1)2^m}\right\},\quad\text{and}\\
&S:=\left\{\left|T\left(\vv{\nu^N}_{1,l},\vv{g}_{l+1,m}\right)\right|>\frac{t}{(l+1)2^m}\right\}.
\end{align*}
We will control each $|S_k|$ and $|S|$ individually.

We will first control $|S_k|$. Begin by using Chebyshev's inequality, the fact that $(b_{k,j_k}^Ndm-v_{k,j_k}^N)(Q_{k,j_k})=0$, Fubini, and trivial estimates to see
\begin{align*}
|S_k|&\leq (l+1)2^mt^{-1}\int_{\mathbb{R}^n\setminus G}\left|T\left(\vv{\nu^N}_{1,k-1},b_{k}^N-\nu_{k}^N,\vv{b^N}_{k+1,l},\vv{g}_{l+1,m}\right)(x)\right|dx\\
&\leq (m+1)2^mt^{-1}\sum_{j_1,\ldots,j_l=1}^N\int_{\mathbb{R}^n\setminus G}\left|\int_{\mathbb{R}^{n(m-l)}}\int_{Q_{l,j_l}}\cdots\int_{Q_{1,j_1}}K(x,\vv{y}_{1,m})\right.\\
&\quad\quad \times \left.\left(\prod_{i=k+1}^lb_{i,j_i}(y_i)\right)\left(\prod_{i=l+1}^mg_i(y_i)\right)\right.\\
&\quad\quad\quad \left. d\vv{\nu}_{(1,j_1),(k-1,j_{k-1})}(\vv{y}_{1,k-1})d(b_{k,j_k}dm-\nu_{k,j_k})(y_k)d\vv{y}_{k+1,m}\right|dx\\
&\leq (m+1)2^mt^{-1}\sum_{j_1,\ldots,j_l=1}^N\int_{\mathbb{R}^{n(m-l)}}\int_{Q_{l,j_l}}\cdots\int_{Q_{1,j_1}}\int_{\mathbb{R}^n\setminus G}\\
&\quad\quad \times |K(x,\vv{y}_{1,m})-K(x,\vv{c}_{(1,j_1),(l,j_l)},\vv{y}_{l+1,m})|\left(\prod_{i=k+1}^l|b_{i,j_i}(y_i)|\right)\left(\prod_{i=l+1}^m|g_i(y_i)|\right)\\
&\quad\quad\quad dxd\vv{|\nu|}_{(1,j_1),(k-1,j_{k-1})}(\vv{y}_{1,k-1})d|b_{k,j_k}dm-\nu_{k,j_k}|(y_k)d\vv{y}_{k+1,m}\\
&\leq (m+1)2^mt^{-1}\sum_{j_1,\ldots,j_l=1}^N\left(\prod_{i=1}^{k-1}|a_{i,j_i}|\right)|b_{k,j_k}dm-\nu_{k,j_k}|(Q_{k,j_k})\\
&\quad\quad\times\left(\prod_{i=k+1}^l\|b_{i,j_i}\|_{L^1(\mathbb{R}^n)}\right)\left(\prod_{i=l+1}^m\|g_i\|_{L^{\infty}(\mathbb{R}^n)}\right)\\
&\quad\quad\times\int_{\mathbb{R}^{n(m-l)}}\sup_{\substack{(y_1,\ldots,y_l) \\\in \prod_{i=1}^lQ_{i,j_i}}} \int_{\mathbb{R}^n\setminus G} |K(x,\vv{y}_{1,m})-K(x,\vv{c}_{(1,j_1),(l,j_l)},\vv{y}_{l+1,m})|dxd\vv{y}_{l+1,m}.
\end{align*}

Now, note that $\|b_{i,j}\|_{L^1(\mathbb{R}^n)}\leq (17\sqrt{n})^n t^{\frac{1}{m}}|Q_{i,j}|$. Indeed, for a fixed $Q_{i,j}$, let $Q_{i,j}^*$ be the cube with the same center, but diameter $17\sqrt{n}$ times as large. Then $Q_{i,j}^*\cap (\mathbb{R}^n\setminus G_i) \neq \emptyset$. So there is a point $x\in Q_{i,j}^*$ such that $Mf_i(x)\leq t^{\frac{1}{m}}$. In particular, $\int_{Q_{i,j}^*}|f_i(y)|dy\leq t^{\frac{1}{m}}|Q_{i,j}^*|$. Since $|Q_{i,j}^*|=(17\sqrt{n})^n|Q_{i,j}|$, we have \[\|b_{i,j}\|_{L^1(\mathbb{R}^n)}=\int_{Q_{i,j}}|f_i(y)|dy\leq\int_{Q_{i,j}^*}|f_i(y)|dy\leq t^{\frac{1}{m}}|Q_{i,j}^*| = (17\sqrt{n})^nt^{\frac{1}{m}}|Q_{i,j}|.\] Use the fact that $|b_{k,j_k}^Ndm-v_{k,j_k}^N|(Q_{k,j_k})\leq 2\|b_{k,j_k}\|_{L^1(\mathbb{R}^n)}$, the $L^{\infty}$ control of the $g_i$, the $L^1$ control of the $b_{i,j}$, and Lemma 1 (which applies since $2\text{diam}(Q_{i,j})\leq d(Q_{i,j},\mathbb{R}^n\setminus G_i)$) to continue the estimate
\begin{align*}
|S_k|&\leq (m+1)2^{m+1}t^{-1}\sum_{j_1,\ldots,j_l=1}^N\left(\prod_{i=1}^l\|b_{i,j_i}\|_{L^1(\mathbb{R}^n)}\right)\left(\prod_{i=l+1}^m\|g_i\|_{L^{\infty}(\mathbb{R}^n)}\right)\\
&\quad\quad\times\int_{\mathbb{R}^{n(m-l)}}\sup_{\substack{(y_1,\ldots,y_l) \\\in \prod_{i=1}^lQ_{i,j_i}}} \int_{\mathbb{R}^n\setminus G} |K(x,\vv{y}_{1,m})-K(x,\vv{c}_{(1,j_1),(l,j_l)},\vv{y}_{l+1,m})|dxd\vv{y}_{l+1,m}\\
&\leq (m+1)2^{m+1}(17\sqrt{n})^{nl}\sum_{j_1,\ldots,j_l=1}^N\left(\prod_{i=1}^l|Q_{i,j_i}|\right)\\
&\quad\quad\times\int_{\mathbb{R}^{n(m-l)}}\sup_{\substack{(y_1,\ldots,y_l) \\\in \prod_{i=1}^lQ_{i,j_i}}} \int_{\mathbb{R}^n\setminus G} |K(x,\vv{y}_{1,m})-K(x,\vv{c}_{(1,j_1),(l,j_l)},\vv{y}_{l+1,m})|dxd\vv{y}_{l+1,m}\\
&\leq A_1(m+1)2^{m+1}(17\sqrt{n})^{nm}\sum_{i=1}^l|G_i|\\
&\leq A_1m(m+1)2^{m+1}(17\sqrt{n})^{nm}\|M\|t^{-\frac{1}{m}}.
\end{align*}

The control of $|S|$ follows from applying Theorem 2 below
\begin{align*}
|S|&\leq 2(l+1)^{\frac{1}{m}}A_3t^{-\frac{1}{m}}\left(\prod_{i=1}^l\left\|\nu_i^N\right\|^{\frac{1}{m}}\right)\left(\prod_{i=l+1}^m\|g_{i}\|_{L^1(\mathbb{R}^n)}^{\frac{1}{m}}\right)\\
&\leq 2(m+1)^{\frac{1}{m}}A_3t^{-\frac{1}{m}}\left(\prod_{i=1}^l\left\|b_i^N\right\|_{L^1(\mathbb{R}^n)}^{\frac{1}{m}}\right)\left(\prod_{i=l+1}^m\|g_{i}\|_{L^1(\mathbb{R}^n)}^{\frac{1}{m}}\right)\\
&\leq 2(m+1)^{\frac{1}{m}}A_3t^{-\frac{1}{m}}\left(\prod_{i=1}^m\|f_i\|_{L^1(\mathbb{R}^n)}^{\frac{1}{m}}\right)\\
&=2(m+1)^{\frac{1}{m}}A_3t^{-\frac{1}{m}}.
\end{align*}
Put the estimates of $|S_k|$ and $|S|$ together to get
\begin{align*}
|\widetilde{E}_s|&\leq \left(m^2\|M\|+\sum_{k=1}^lA_1m(m+1)2^{m+1}(17\sqrt{n})^{nm}\|M\|+2(m+1)^{\frac{1}{m}}A_3\right)t^{-\frac{1}{m}}\\
&\leq \left(m^2\|M\|+A_1m^2(m+1)2^{m+1}(17\sqrt{n})^{nm}\|M\|+2(m+1)^{\frac{1}{m}}A_3\right)t^{-\frac{1}{m}}\\
&=B_2t^{-\frac{1}{m}}
\end{align*}
where $B_2:=m^2\|M\|+A_1m^2(m+1)2^{m+1}(17\sqrt{n})^{nm}\|M\|+2(m+1)^{\frac{1}{m}}A_3$. Since the above estimate is independent of $N$, letting $N\rightarrow \infty$ yields \[|E_s| \leq B_2t^{-\frac{1}{m}}\]

Finally, use the estimates of $|E_s|$, $1\leq s \leq 2^m$ to observe 
\[|\{|T(f_1,\ldots,f_m)|>t\}| \leq |E_1|+\sum_{s=2}^{2^m}|E_s|\leq \left(B_1+\left(2^m-1\right)B_2\right)t^{-\frac{1}{m}}.\]
Take $A_2=\left(B_1+\left(2^m-1\right)B_2\right)^m$ to complete the proof.
\end{proof}

We now prove Theorem 2.

\begin{proof}[Proof of Theorem 2]
Assume without loss of generality that each $a_{i,j}>0$ and that $\|\nu_1\|=\cdots=\|\nu_l\|=\|f_{l+1}\|_{L^1(\mathbb{R}^n)}=\cdots=\|f_{m}\|_{L^1(\mathbb{R}^n)}=1$. For $i=1,\ldots,l$, set $$E_{i,1}:=B(x_{i,1},r_{i,1})$$ where $r_{i,1}>0$ is chosen so that $|E_{i,1}|=a_{i,1}t^{-\frac{1}{m}}$. Subsequently, set $$E_{i,2}:= B(x_{i,2},r_{i,2})\setminus E_{i,1}$$ where $r_{i,2}>0$ is chosen so that $|E_{i,2}|=a_{i,2}t^{-\frac{1}{m}}$. In general, for $j=3,\ldots, N$, set $$E_{i,j}:=B(x_{i,j},r_{i,j})\setminus \bigcup_{k=1}^{j-1}E_{i,k}$$ where $r_{i,j}>0$ is chosen so that $|E_{i,j}|=a_{i,j}t^{-\frac{1}{m}}$. Set $$E_i:=\bigcup_{j=1}^{N}E_{i,j}$$ and notice that by construction $$|E_i|=\sum_{j=1}^N|E_{i,j}|=\sum_{j=1}^Na_{i,j}t^{-\frac{1}{m}}=\|\nu_i\|t^{-\frac{1}{m}}=t^{-\frac{1}{m}}.$$ Similarly, set $$E_{i,1}^* := B(x_{i,1},2r_{i,1}),$$ and subsequently for $j=2,\ldots,N$, $$E_{i,j}^*:=B(x_{i,j},2r_{i,j})\setminus \bigcup_{k=1}^{j-1}E_{i,k}^*.$$ Set $$\quad E_i^*:=\bigcup_{j=1}^NE_{i,j}^*, \quad \text{and} \quad E^*:=\bigcup_{i=1}^{l}E_{i}^*.$$ By the doubling property of Lebesgue measure, $$|E^*|\leq \sum_{i=1}^l|E^*_{i}|\leq 2^n\sum_{i=1}^l|E_{i}|\leq m2^{n}t^{-\frac{1}{m}}.$$

For $k \in \{0,\ldots,l\}$, set 
\[\sigma_k :=T\left(t^{\frac{1}{m}}\mathbbm{1}_{E_{1}},\ldots,t^{\frac{1}{m}}\mathbbm{1}_{E_{k}},\nu_{k+1}\ldots,\nu_l,f_{l+1},\ldots,f_{m}\right),\] 
noticing that $\sigma_0=T(\nu_1,\ldots,\nu_{l},f_{l+1},\ldots,f_{m})$.
Then, by adding and subtracting $\sigma_k$ for $1\leq k\leq l$, we have
\begin{align*}
|\{|T(\nu_1,\ldots&,\nu_{l},f_{l+1},\ldots,f_{m})|>t\}|\leq\sum_{k=1}^{l}\left|\left\{|\sigma_{k-1}-\sigma_{k}|>\frac{t}{l+1}\right\}\right|+\left|\left\{\left|\sigma_l\right|>\frac{t}{l+1}\right\}\right|\\
&=\sum_{k=1}^{l}|E^*|+\sum_{k=1}^{l}\left|\left\{\mathbb{R}^n\setminus E^*:|\sigma_{k-1}-\sigma_{k}|>\frac{t}{l+1}\right\}\right|+\left|\left\{\left|\sigma_l\right|>\frac{t}{l+1}\right\}\right|\\
&\leq m^22^nt^{-\frac{1}{m}}+\sum_{k=1}^{l}|P_k|+|P|,
\end{align*}
where 
\begin{align*}
&P_k:=\left\{\mathbb{R}^n\setminus E^*:|\sigma_{k-1}-\sigma_{k}|>\frac{t}{l+1}\right\},\quad\quad\text{and}\\
&P:=\left\{\left|\sigma_l\right|>\frac{t}{l+1}\right\}.
\end{align*}
The remainder of the proof will focus on bounding $|P|$ and each $|P_k|$ by a constant multiplied by $t^{-\frac{1}{m}}$.

To control $|P|$, use Chebyshev's inequality, the boundedness of $T$ from $(L^2(\mathbb{R}^n))^m$ to $L^{\frac{2}{m}}(\mathbb{R}^n)$, and the $L^{\infty}$ control of the $f_i$ to observe
\begin{align*}
|P|&\leq \frac{(l+1)^{\frac{2}{m}}}{t^{\frac{2}{m}-\frac{2l}{m^2}}}\int_{\mathbb{R}^n}\left|T\left(\mathbbm{1}_{E_1},\ldots,\mathbbm{1}_{E_l},f_{l+1},\ldots,f_{m}\right)(x)\right|^{\frac{2}{m}}dx\\
&\leq \frac{(m+1)^{\frac{2}{m}}}{t^{\frac{2}{m}-\frac{2l}{m^2}}}\|T\|^{\frac{2}{m}}_{(L^2(\mathbb{R}^n))^m\rightarrow L^{\frac{2}{m}}(\mathbb{R}^n)}\left(\prod_{i=1}^l\|\mathbbm{1}_{E_i}\|_{L^2(\mathbb{R}^n)}^{\frac{2}{m}}\right)\left(\prod_{i=l+1}^{m}\|f_i\|_{L^2(\mathbb{R}^n)}^{\frac{2}{m}}\right)\\
&\leq \frac{(m+1)^{\frac{2}{m}}}{t^{\frac{2}{m}-\frac{2l}{m^2}}}\|T\|^{\frac{2}{m}}_{(L^2(\mathbb{R}^n))^m\rightarrow L^{\frac{2}{m}}(\mathbb{R}^n)}\left(t^{-\frac{l}{m^2}}\right)\left(t^{\frac{m-l}{m^2}}\right)\prod_{i=l+1}^{m}\|f_i\|_{L^1(\mathbb{R}^n)}^{\frac{1}{m}}\\
&\leq (m+1)^{\frac{2}{m}}\|T\|^{\frac{2}{m}}_{(L^2(\mathbb{R}^n))^m\rightarrow L^{\frac{2}{m}}(\mathbb{R}^n)}t^{-\frac{1}{m}}.
\end{align*}

We will now control $|P_k|$. Notice
\[\sigma_{k-1}-\sigma_k=T\left(t^{\frac{1}{m}}\mathbbm{1}_{E_1},\ldots,t^{\frac{1}{m}}\mathbbm{1}_{E_{k-1}},\nu_{k}-t^{\frac{1}{m}}\mathbbm{1}_{E_k},\nu_{k+1},\ldots,\nu_l,f_{l+1},\ldots,f_m\right).\]
Use Chebyshev's inequality, the fact that $(\nu_{k,j_k}-t^{\frac{1}{m}}\mathbbm{1}_{E_{k,j_k}}dm)(E_k)=0$, Fubini, and trivial estimates to see
\begin{align*}
&|P_k|\leq \frac{l+1}{t^{\frac{m-k+1}{m}}}\int_{\mathbb{R}^n\setminus E^*}\left|T\left(\mathbbm{1}_{E_{1}},\ldots,\mathbbm{1}_{E_{k-1}},\nu_k-t^{\frac{1}{m}}\mathbbm{1}_{E_k},\nu_{k+1},\ldots,\nu_l,f_{l+1},\ldots,f_m\right)(x)\right|dx\\
&\leq \frac{m+1}{t^{\frac{m-k+1}{m}}}\sum_{j_1,\ldots,j_l=1}^N\int_{\mathbb{R}^n\setminus E^*}\left|\int_{\mathbb{R}^{n(m-l)}}\int_{E_{l,j_l}}\cdots\int_{E_{1,j_1}}K(x,\vv{y}_{1,m})\right.\\
&\quad\quad\times\left.\left(\prod_{i=l+1}^mf_{i}(y_{i})\right)d\vv{y}_{1,k-1}d(\nu_{k,j_k}-t^{\frac{1}{m}}\mathbbm{1}_{E_{k,j_k}}dm)(y_k)d\vv{\nu}_{(k+1,j_{k+1}),(l,j_l)}(\vv{y}_{k+1,l})d\vv{y}_{l+1,m}\right|dx\\
&\leq \frac{m+1}{t^{\frac{m-k+1}{m}}}\sum_{j_1,\ldots,j_l=1}^N\int_{\mathbb{R}^{n(m-l)}}\int_{E_{l,j_l}}\cdots\int_{E_{1,j_1}}\int_{\mathbb{R}^n\setminus E^*}|K(x,\vv{y}_{1,m})-K(x,\vv{x}_{(1,j_1),(l,j_l)},\vv{y}_{l+1,m})|\\
&\quad\quad\times\left(\prod_{i=l+1}^m|f_{i}(y_{i})|\right)dxd\vv{y}_{1,k-1}d\left|\nu_{k,j_k}-t^{\frac{1}{m}}\mathbbm{1}_{E_{k,j_k}}dm\right|(y_k)d\vv{|\nu|}_{(k+1,j_{k+1}),(l,j_l)}(\vv{y}_{k+1,l})d\vv{y}_{l+1,m}\\
&\leq \frac{m+1}{t^{\frac{m-k+1}{m}}}\sum_{j_1,\ldots,j_l=1}^N\left(\prod_{i=1}^{k-1}|E_{i,j_i}|\right)\left|\nu_{k,j_k}-t^{\frac{1}{m}}\mathbbm{1}_{E_{k,j_k}}dm\right|(E_{k,j_k})\left(\prod_{i=k+1}^la_{i,j_i}\right)\left(\prod_{i=l+1}^m\|f_i\|_{L^{\infty}(\mathbb{R}^n)}\right)\\
&\quad\quad\times\int_{\mathbb{R}^{n(m-l)}}\sup_{\substack{(y_1,\ldots,y_l) \\\in \prod_{i=1}^lE_{i,j_i}}}\int_{\mathbb{R}^n\setminus E^*}|K(x,\vv{y}_{1,m})-K(x,\vv{x}_{(1,j_1),(l,j_l)},\vv{y}_{l+1,m})|dxd\vv{y}_{l+1,m}.
\end{align*}
Use the fact that $|\nu_{k,j_k}-t^{\frac{1}{m}}\mathbbm{1}_{E_{k,j_k}}dm|(E_k)\leq 2t^{\frac{1}{m}}|E_{k,j_k}|$, the $L^{\infty}(\mathbb{R}^n)$ control of the $f_i$, and Lemma 1 to continue the estimate 
\begin{align*}
&|P_k| \leq 2(m+1)\sum_{j_1,\ldots,j_l=1}^N\left(\prod_{i=1}^{l}|E_{i,j_i}|\right)\\
&\quad\quad\times\int_{\mathbb{R}^{n(m-l)}}\sup_{\substack{(y_1,\ldots,y_l) \\\in \prod_{i=1}^lE_{i,j_i}}}\int_{\mathbb{R}^n\setminus E^*}|K(x,\vv{y}_{1,m})-K(x,\vv{x}_{(1,j_1),(l,j_l)},\vv{y}_{l+1,m})|dxd\vv{y}_{l+1,m}\\
&\leq 2(m+1)A_1\sum_{i=1}^l|E_i|\\
&\leq 2m(m+1)A_1t^{-\frac{1}{m}}.
\end{align*}

Using these estimates of $|P|$ and $|P_k|$, we have 
\begin{align*}
|\{|T(\nu_1,\ldots,\nu_l&,f_1,\ldots,f_{m-l})|>t\}|\leq m^22^nt^{-\frac{1}{m}}+\sum_{k=1}^{l}|P_k|+|P|\\
&\leq\left(m^22^n+\left(\sum_{k=1}^{l}2m(m+1)A_1\right)+(m+1)^{\frac{2}{m}}\|T\|^{\frac{2}{m}}_{(L^2(\mathbb{R}^n))^m\rightarrow L^{\frac{2}{m}}(\mathbb{R}^n)}\right)t^{-\frac{1}{m}}\\
&\leq\left(m^22^n+2m^2(m+1)A_1+(m+1)^{\frac{2}{m}}\|T\|^{\frac{2}{m}}_{(L^2(\mathbb{R}^n))^m\rightarrow L^{\frac{2}{m}}(\mathbb{R}^n)}\right)t^{-\frac{1}{m}}.
\end{align*}
Take $A_3=\left(m^22^n+2m^2(m+1)A_1+(m+1)^{\frac{2}{m}}\|T\|^{\frac{2}{m}}_{(L^2(\mathbb{R}^n))^m\rightarrow L^{\frac{2}{m}}(\mathbb{R}^n)}\right)$ to complete the proof.
\end{proof}


\end{document}